\documentclass{amsart}

\usepackage{amsmath,amsfonts,amssymb}
\usepackage{tikz}
\usepackage{hyperref}

\usepackage{xcolor}
\usepackage[colorinlistoftodos]{todonotes}

\newcommand{\svw}[1]{\textcolor{black}{#1}}

\newtheorem{theorem}{Theorem}[section]
\newtheorem{lemma}[theorem]{Lemma}
\newtheorem{proposition}[theorem]{Proposition}
\newtheorem{corollary}[theorem]{Corollary}

\theoremstyle{definition}
\newtheorem{definition}[theorem]{Definition}
\newtheorem{example}[theorem]{Example}

\theoremstyle{remark}
\newtheorem{remark}[theorem]{Remark}

\numberwithin{equation}{section}

\newcommand{\newword}[1]{\emph{\textbf{#1}}}

\newcommand{\key}{\kappa}
\newcommand{\fund}{\mathfrak{F}}

\newcommand{\SKT}{\mathrm{SKT}}

\newcommand{\des}{\mathrm{des}}
\newcommand{\comp}[1]{\mathbf{#1}}

\newcommand{\skews}{/}

\newcommand{\rk}{\mathrm{rk}}

\newcommand{\cover}{\,{\buildrel \over \prec\!\!\!\cdot}\,}

\newlength\cellsize 

\savebox2{%
\begin{picture}(10,10)
\put(0,0){\line(1,0){10}}
\put(0,0){\line(0,1){10}}
\put(10,0){\line(0,1){10}}
\put(0,10){\line(1,0){10}}
\end{picture}}

\newcommand\cellify[1]{\def\thearg{#1}\def\nothing{}%
\ifx\thearg\nothing\vrule width0pt height\cellsize depth0pt%
  \else\hbox to 0pt{\usebox2\hss}\fi%
  \vbox to \cellsize{\vss\hbox to \cellsize{\hss$_{#1}$\hss}\vss}}

\newcommand\tableau[1]{\vtop{\let\\=\cr
\setlength\baselineskip{-10000pt}
\setlength\lineskiplimit{10000pt}
\setlength\lineskip{0pt}
\halign{&\cellify{##}\cr#1\crcr}}}

\newlength\smcellsize 

\savebox3{%
\begin{picture}(8,8)
\put(0,0){\line(1,0){8}}
\put(0,0){\line(0,1){8}}
\put(8,0){\line(0,1){8}}
\put(0,8){\line(1,0){8}}
\end{picture}}

\newcommand\smcellify[1]{\def\thearg{#1}\def\nothing{}%
\ifx\thearg\nothing\vrule width0pt height\smcellsize depth0pt%
  \else\hbox to 0pt{\usebox3\hss}\fi%
  \vbox to \smcellsize{\vss\hbox to \smcellsize{\hss$_{#1}$\hss}\vss}}

\newcommand\smtab[1]{\vtop{\let\\=\cr
\setlength\baselineskip{-8000pt}
\setlength\lineskiplimit{8000pt}
\setlength\lineskip{0pt}
\halign{&\smcellify{##}\cr#1\crcr}}}

\newcommand\cbox[1]{\color{#1!30}\rule{1\smcellsize}{1\smcellsize}\hspace{-\smcellsize}\usebox3}

\begin{document}


\title[Skew key polynomials]{Skew key polynomials and a generalized Littlewood--Richardson rule}  

\author[S. Assaf]{Sami Assaf}
\address{Department of Mathematics, University of Southern California,  3620 S. Vermont Ave., Los Angeles, CA 90089-2532, U.S.A.}
\email{shassaf@usc.edu}
\thanks{S.A. supported in part by NSF DMS-1763336. S.v.W. supported in part by NSERC}

\author[S. van Willigenburg]{Stephanie van Willigenburg}
\address{Department of Mathematics, University of British Columbia, 1984 Mathematics Rd., Vancouver, BC V6T 1Z2, Canada}
\email{steph@math.ubc.ca}





\begin{abstract}
  Young's lattice is a partial order on integer partitions whose saturated chains correspond to standard Young tableaux, \svw{one type of combinatorial object that generates} the Schur basis for symmetric functions. Generalizing Young's lattice, we introduce a new partial order on weak compositions that we call the key poset. Saturated chains in this poset correspond to standard key tableaux, the combinatorial objects that generate the key polynomials, a nonsymmetric polynomial generalization of the Schur basis. Generalizing skew Schur functions, we define skew key polynomials in terms of this new poset. Using weak dual equivalence, we give a nonnegative weak composition Littlewood--Richardson rule for the key expansion of skew key polynomials, generalizing the flagged Littlewood--Richardson rule of Reiner and Shimozono.
\end{abstract}

\maketitle

%
\section{Introduction}
%
\label{sec:introduction}

Schur polynomials are central to the study of the representation theories of the general linear group and of the symmetric group, as well as to understanding the geometry of the Grassmannian. The combinatorics of Young tableaux, the ubiquitous objects that generate the Schur polynomials, often sheds light on important representation theoretic or geometric properties such as tensor products, induction and restriction of modules, and intersection multiplicities.

The celebrated Littlewood--Richardson rule \cite{LR34} gives a combinatorial description for the Schur expansion of a product of two Schur polynomials or, equivalently, for the Schur expansion of a \emph{skew Schur polynomial} as
\[ s_{\lambda} s_{\mu} = \sum_{\nu} c_{\lambda,\mu}^{\nu} s_{\nu},
\hspace{1em} \text{or} \hspace{1em}
s_{\nu\skews\lambda} = \sum_{\mu} c_{\lambda,\mu}^{\nu} s_{\mu}, \]
where $c_{\lambda,\mu}^{\nu}$ is the number of saturated chains in \emph{Young's lattice} from $\lambda$ to $\nu$ satisfying certain conditions depending on $\mu$. Here Young's lattice is the partial order on integer partitions given by containment of Young diagrams. These so-called \emph{Littlewood--Richardson coefficients} arise in representation theory as the irreducible multiplicities for the tensor product of two irreducible representations for the general linear group and as the irreducible multiplicities for the induced tensor product of two irreducible representations for the symmetric group. They also appear geometrically, giving the number of points lying in a suitable intersection of three Grassmannian Schubert varieties.

The \emph{key polynomials} are nonsymmetric polynomial generalizations of Schur polynomials first studied by Demazure \cite{Dem74a} in connection with Schubert varieties. Key polynomials are irreducible characters of Demazure modules for the general linear group \cite{Dem74} and represent Schubert classes for vexillary permutations \cite{LS90}. They form an important basis for the polynomial ring, and so we may consider their structure constants parallel to those for the Schur basis of symmetric polynomials. However, in stark contrast with the Schur case, the structure constants of key polynomials are not, in general, nonnegative, though Kouno \cite{Kuono} has partial results where nonnegativity holds and Assaf and Quijada \cite{AQ18} have made progress on understanding the signs in the Pieri case.

Nevertheless, key polynomials appear in many of the myriad generalizations of the Littlewood--Richardson rule. One such rule is the \emph{flagged Littlewood--Richardson rule} of Reiner and Shimozono \cite{RS95}. Flagged Schur polynomials \cite{LS82} arise as those polynomials occuring both as key polynomials and as Schubert polynomials. Reiner and Shimozono \cite{RS95} considered the \emph{flagged skew Schur polynomials} and gave a nonnegative rule for their expansion into the key basis. Assaf \cite{Ass-W} considered a diagram containment-based skew analog of key polynomials along the same lines, but obtained nonnegativity results only in very special cases.

Another nonnegative rule is the \emph{quasisymmetric Littlewood--Richardson rule}. The quasisymmetric Schur functions of Haglund, Luoto, Mason and van Willigenburg \cite{HLMvW11-2} are a quasisymmetric generalization of Schur polynomials whose combinatorics shares many nice properties with that for Schur polynomials. These authors \cite{HLMvW11} use the quasisymmetric Schur functions to derive a nonnegative \emph{refined Littlewood--Richardson rule} for the product of a \emph{key polynomial} and a Schur polynomial with sufficiently many variables. Related to this, Bessenrodt, Luoto and van Willigenburg \cite{BLvW11} define a \emph{partial order} on strong compositions that gives rise to a nonnegative Littlewood--Richardson rule for a skew analog of the quasisymmetric Schur functions. For details on these results and quasisymmetric Schur functions in general, see the book by Luoto, Mykytiuk, and van Willigenburg \cite{LMvW}.

In this paper, we generalize the flagged Littlewood--Richardson rule to skew key polynomials in the largest possible setting where nonnegativity prevails. To do so, we begin in Section~\ref{sec:poset} by generalizing Young's lattice to a partial order on \emph{weak} compositions that we call the \emph{key poset}. In contrast with the generalization to strong compositions in \cite{BLvW11}, we give explicit cover relations as well as explicit criteria for comparability in the poset, 
\svw{and remark on the latter hence filling this gap}, 
though as with the strong composition poset, the key poset is not a lattice. In Section~\ref{sec:poly}, we relate the poset with the tableaux combinatorics for key polynomials. Using these {paradigms} together, in Section~\ref{sec:skew} we re-define skew key polynomials with respect to the poset and give a general nonnegative Littlewood--Richardson rule for skew key polynomials, vastly generalizing the nonnegativity results of \cite{Ass-W, RS95}.
The relations in the key poset are more \svw{restrictive} than mere containment used to define skew key polynomials considered in \cite{Ass-W}. Under the more general containment \svw{definition \cite{Ass-W},} all key polynomial coefficients are nonnegative if and only if the two indexing shapes are comparable in the key poset. Thus the key poset is precisely the right notion to capture positivity.

%
\section{Posets}
%
\label{sec:poset}


An integer \newword{partition} $\lambda = (\lambda_1,\lambda_2,\ldots,\lambda_{\ell})$ is a weakly decreasing sequence of positive integers, $\lambda_i \geq \lambda_{i+1} > 0$. The \newword{rank} of a partition $\lambda$, denoted by $\rk(\lambda)$, is the sum of the parts, 
\[ \rk(\lambda) = \lambda_1 + \lambda_2 + \cdots + \lambda_{\ell}  \]{and we call $\ell$ its \newword{length}.}

The \newword{Young diagram} of a partition $\lambda$ is the collection of $\lambda_i$ unit cells left-justified in row $i$ indexed from the bottom (French notation). Abusing notation, we use $\lambda$ interchangeably for the integer partition and for its diagram.

A \newword{partially ordered set}, or \newword{poset}, is a set together with a partial order comparing certain elements of the set. We turn integer partitions into a poset by the containment relation, setting $\lambda \subseteq \mu$ if $\lambda_i \leq \mu_i$ for all $i$ or, equivalently, if the diagram for $\lambda$ is a subset of the diagram for $\mu$. We call this poset \newword{Young's lattice}.

For $p,q$ in a poset $\mathcal{P}$, we say $q$ \newword{covers} $p$, denoted by $p \cover q$, if $p \prec q$ and for any $r\in\mathcal{P}$ for which $p \preceq r \preceq q$, either $p=r$ or $r=q$. The \newword{cover relations} for Young's lattice may be described by $\lambda \cover \mu$ if and only if $\mu$ is obtained from $\lambda$ by incrementing a single part {$\lambda_{i+1}$} for which $\lambda_i > \lambda_{i+1}$ by $1$ or, equivalently, by adding a single {cell} to the end of a row for which the {row with one smaller index is strictly longer.}

A poset $\mathcal{P}$ is a \newword{lattice} if every pair of elements $p,q\in\mathcal{P}$ have a unique least upper bound and a unique greatest lower bound. For Young's lattice, these constructions are given by the set-theoretic union and the set-theoretic intersection of the diagrams, respectively. 

Young's lattice is a prominent tool in algebraic combinatorics, used to study symmetric functions, representations of finite and affine Lie groups, and intersection numbers for finite and affine Grassmannians. We generalize the construction from integer partitions to weak compositions in such a way that maintains the connection to representation theory and geometry.


A \newword{weak composition} $\comp{a} = (a_1,a_2,\ldots,a_{n})$ is a sequence of nonnegative integers, that is, $a_i \geq 0$. {A weak composition is a \newword{strong composition} if all $a_i>0$. Given a weak composition $\comp{a}$ we denote the strong composition obtained by removing its zeros by $\mathrm{flat}(\comp{a})$.} Extending notation, the \newword{rank} of a weak composition $\comp{a}$, denoted by $\rk(\comp{a})$, is the sum of the parts, 
\[ \rk(\comp{a}) = a_1 + a_2 + \cdots + a_{n}  \]{and we call $n$ its \newword{length}.}

The \newword{key diagram} of a weak composition $\comp{a}$ is the collection of $a_i$ unit cells left-justified in row $i$ indexed from the bottom. As above, we use $\comp{a}$ interchangeably for the weak composition and for its {key} diagram.

\begin{definition}
  The \newword{key poset} is the partial order $\prec$ on weak compositions of length $n$ defined by the relation $\comp{a} \preceq \comp{b}$ if and only if $a_i \leq b_i$ for $i=1,2,\ldots,n$ and for any indices $1 \leq i<j \leq n$ for which $b_j > a_j$ and $a_i > a_j$, we have $b_i > b_j$.
  \label{def:above}
\end{definition}

{If $\comp{a} \preceq \comp{b}$, then the collection of cells in $\comp{b}$ but not in $\comp{a}$ is called a \newword{skew key diagram} denoted by $\comp{b}\skews \comp{a}$.}  We will also consider more general skew diagrams whenever $\comp{a} \subseteq \comp{b}$, or $a_i \leq b_i$ for all $i$.

{Definition~\ref{def:above} has the following interpretation in terms of key diagrams that we will use often:}  $\comp{a} \prec \comp{b}$ if and only if $\comp{a} \subset \comp{b}$ and whenever a cell of $\comp{b}\skews \comp{a}$ lies above a cell of $\comp{a}$, the lower row is strictly longer in $\comp{b}$; see Fig.~\ref{fig:above}. 

\begin{figure}[ht]
  \begin{displaymath}
    \vline\smtab{ %
      \cbox{white} & \cbox{white} & \cbox{blue} & \cbox{blue} \\ \\
      \cbox{white} & \cbox{white} & \cbox{white} & \cbox{purple} & \cbox{purple}}
  \end{displaymath}
  \caption{\label{fig:above}An illustration of the partial order on weak compositions in terms of key diagrams. Here, cells $\cbox{white}$ lie in $\comp{a} \subseteq \comp{b}$, cells $\cbox{blue}$ lie in $\comp{b} \skews \comp{a}$, and cells $\cbox{purple}$ lie in $\comp{b} = \comp{a} \cup \comp{b}$.}
\end{figure}

Notice this partial order is not given simply by containment of key diagrams.

\begin{example}
  Let $\comp{a}=(2,1,1)$ and $\comp{b}=(2,1,2)$. Then $a_i \leq b_i$ for $i=1,2,3$, showing $\comp{a} \subset \comp{b}$. However, $\comp{a} \not\prec \comp{b}$ since $b_3 > a_3$ and $a_1 > a_3$ but $b_1 = b_3$. See Fig.~\ref{fig:contain}.
  \label{ex:contain}
\end{example}

\begin{figure}[ht]
  \begin{displaymath}
    \vline\smtab{ \ \\ \ \\ \ & \ }
    \hspace{3\cellsize} 
    \vline\smtab{ \ & \cbox{purple} \\ \ \\ \ & \ }
  \end{displaymath}
  \caption{\label{fig:contain}Two incomparable diagrams where the left is contained in the right as a subset but is not below in the poset as evidenced by the set-theoretic difference shown as $\cbox{purple}$.}
\end{figure}

As with Young's lattice, this partial order on weak compositions is ranked by the number of {cells}, and we may describe the covering relations in terms of adding a cell subject to certain conditions.

\begin{theorem}
  The {key} poset on weak compositions is ranked by $\rk$ with covering relation $\comp{a} \cover \comp{b}$ if and only if $\comp{b}$ is obtained from $\comp{a}$ by incrementing $a_j$ by $1$ where for any $i<j$ we have $a_i \neq a_j + 1$.
  \label{thm:cover}
\end{theorem}

\begin{proof}
  Let $\comp{a} \prec \comp{b}$ be two comparable weak compositions of length $n$. If $\rk(\comp{b}) = \rk(\comp{a})+1$, then $\comp{a}\subset \comp{b}$ and letting $j$ denote the row index of the unique cell of $\comp{b} \skews \comp{a}$, then $b_j>a_j$ and for any $i<j$ for which $a_i\geq a_j+1$, we must have $a_i = b_i > b_j = a_j+1$. In particular, we must have either $a_i \leq a_j$ or $a_i > a_j+1$, satisfying the stated cover relation.

  Conversely, given $\comp{a}$ and a row index $j$ for which $a_i \leq a_j$ or $a_i > a_j+1$ for all $i<j$, the weak composition $\comp{b}$ defined by $b_i=a_i$ for $i\neq j$ and $b_j = a_j+1$ satisfies the condition that for any row index $i<j$ for which $a_i > a_j$, we have $b_i = a_i > a_j +1 = b_j$. Thus the cover relation implies $\comp{a} \prec \comp{b}$.

  Finally, suppose $\comp{a} \prec \comp{c}$ with $\rk(\comp{c}) > \rk(\comp{a})+1$. Let $k$ denote the shortest, then lowest if tied, row of $\comp{c}$ containing an element of $\comp{c} \skews \comp{a}$. Set $b_i=c_i$ for $i\neq k$ and $b_k = c_k-1$. Then clearly $\comp{a} \subset \comp{b} \subset \comp{c}$. We claim $\comp{a} \prec \comp{b} \prec \comp{c}$. Given row indices $i<j$, if $b_j>a_j$ and $a_i>a_j$, then we also have $c_j \geq b_j > a_j$, and so since $\comp{a} \prec \comp{c}$, we have $c_i > c_j$. This ensures $k\neq i$ by choice of $k$ as the shortest row of $\comp{c}\skews \comp{a}$ for $k$, and so $b_i = c_i > c_j \geq b_j$ showing $\comp{a} \prec \comp{b}$. If $c_j>b_j$ and $b_i>b_j$, then $j=k$ and {$c_i > c_k=c_j$} by choice of the shortest then lowest row for $k$, and so $\comp{b} \prec \comp{c}$. In particular, $\comp{a}$ is not covered by $\comp{c}$, proving the poset is ranked with cover relations as stated.
\end{proof}

{\begin{remark}\label{rem:compposet} The composition poset \cite[Definition 2.3]{BLvW11} on strong compositions can be described as being ranked by $\rk$ with covering relation $\comp{a} \cover \comp{b}$ if and only if $\comp{b}$ is obtained from $\comp{a}$ by incrementing $a_j$ by $1$ where for any $i<j$ we have $a_i \neq a_j$, \svw{or appending a 1 at the front. Therefore using a proof similar to the first two paragraphs of the above proof we have a description of it analogous to Definition~\ref{def:above}:  The composition poset is the partial order $\prec$ on strong compositions of length $n$ (prepending zeros if lengths differ) defined by the relation $\comp{a} \preceq \comp{b}$ if and only if $a_i \leq b_i$ for $i=1,2,\ldots,n$ and for any indices $1 \leq i<j \leq n$ for which $b_j > a_j$ and $a_i > a_j-1$, we have $b_i > b_j-1$.}
\end{remark}}

In terms of key diagrams, $\comp{a} \cover \comp{b}$ if and only if $\comp{a} \subset \comp{b}$, there is a single cell of $\comp{b}\skews \comp{a}$, and this cell does not sit above any cell lying at the end of its row. Fig.~\ref{fig:Key} depicts the key poset up to rank $3$.

\begin{figure}[ht]
  \begin{center}
    \begin{tikzpicture}[xscale=1.25,yscale=2]
      \node at (4.5,0.25) (K000) {$\varnothing$};
      \node at (3,1) (K001) {$\vline\smtab{\ \\ \\ \\\hline }$};
      \node at (4.5,1) (K010) {$\vline\smtab{\\ \ \\ \\\hline }$};
      \node at (6,1) (K100) {$\vline\smtab{\\ \\ \ \\\hline }$};
      \node at (1.25,2) (K002) {$\vline\smtab{\ & \ \\ \\ \\\hline }$};
      \node at (2.5,2) (K011) {$\vline\smtab{\ \\ \ \\ & \\\hline}$};
      \node at (3.75,2) (K101) {$\vline\smtab{\ \\ \\ \ & \\\hline}$};
      \node at (5.25,2) (K020) {$\vline\smtab{\\ \ & \ \\ & \\\hline}$};
      \node at (6.5,2) (K110) {$\vline\smtab{\\ \ \\ \ & \\\hline}$};
      \node at (7.75,2) (K200) {$\vline\smtab{\\ \\ \ & \ \\\hline}$};
      \node at (0,3) (K003) {$\vline\smtab{\ & \ & \ \\ \\ & \\\hline}$};
      \node at (1,3) (K012) {$\vline\smtab{\ & \ \\ \ \\ & & \\\hline}$};
      \node at (2,3) (K102) {$\vline\smtab{\ & \ \\ \\ \ & & \\\hline}$};
      \node at (3,3) (K021) {$\vline\smtab{\ \\ \ & \ \\ & & \\\hline}$};
      \node at (4,3) (K111) {$\vline\smtab{\ \\ \ \\ \ & & \\\hline}$};
      \node at (5,3) (K201) {$\vline\smtab{\ \\ \\ \ & \ & \\\hline}$};
      \node at (6,3) (K030) {$\vline\smtab{\\ \ & \ & \ \\ & \\\hline}$};
      \node at (7,3) (K120) {$\vline\smtab{\\ \ & \ & \\ \ \\\hline}$};
      \node at (8,3) (K210) {$\vline\smtab{\\ \ \\ \ & \ & \\\hline}$};
      \node at (9,3) (K300) {$\vline\smtab{\\ \\ \ & \ & \ \\\hline}$};
      \draw[thin] (K000.north) -- (K100.south) ;
      \draw[thin] (K000.north) -- (K010.south) ;
      \draw[thin] (K000.north) -- (K001.south) ;
      \draw[thin] (K001.north) -- (K002.south) ;
      \draw[thin] (K001.north) -- (K011.south) ;
      \draw[thin] (K001.north) -- (K101.south) ;
      \draw[thin] (K010.north) -- (K110.south) ;
      \draw[thin] (K010.north) -- (K020.south) ;
      \draw[thin] (K100.north) -- (K200.south) ;
      \draw[thin] (K002.north) -- (K003.south) ;
      \draw[thin] (K002.north) -- (K012.south) ;
      \draw[thin] (K002.north) -- (K102.south) ;
      \draw[thin] (K011.north) -- (K012.south) ;
      \draw[thin] (K011.north) -- (K111.south) ;
      \draw[thin] (K011.north) -- (K021.south) ;
      \draw[thin] (K101.north) -- (K102.south) ;
      \draw[thin] (K101.north) -- (K201.south) ;
      \draw[thin] (K110.north) -- (K120.south) ;
      \draw[thin] (K110.north) -- (K210.south) ;
      \draw[thin] (K020.north) -- (K021.south) ;
      \draw[thin] (K020.north) -- (K120.south) ;
      \draw[thin] (K020.north) -- (K030.south) ;
      \draw[thin] (K200.north) -- (K201.south) ;
      \draw[thin] (K200.north) -- (K210.south) ;
      \draw[thin] (K200.north) -- (K300.south) ;
    \end{tikzpicture}
  \end{center}
  \caption{\label{fig:Key}The Hasse diagram of the key poset up to rank $3$.}
\end{figure}

While containment is not sufficient for covering in general, it is for the partition case. In this way, the key poset generalizes Young's lattice.

\begin{proposition}
  For a weakly increasing weak composition $\comp{a}$, and for $\comp{b}$ any weak composition, we have $\comp{a} \subseteq \comp{b}$ if and only if $\comp{a} \preceq \comp{b}$.
  \label{prop:keypart}
\end{proposition}

\begin{proof}
  For $\comp{a}$ weakly increasing, we never have $a_i>a_j$ for $j>i$, making the latter condition of Definition~\ref{def:above} vacuously true. Thus containment is comparability.
\end{proof}

\begin{corollary}
  Any finite subposet of Young's lattice is a subposet of the key poset.
\end{corollary}

\begin{proof}
  For $\lambda$ a partition of length $\ell\leq n$, let $\comp{a}_{\lambda} = (0,\ldots,0,\lambda_{\ell},\ldots,\lambda_1)$ be the weakly increasing weak composition of length $n$ whose nonzero parts rearrange to $\lambda$. Given any weak composition $\comp{b}$ of length $n$, by Proposition~\ref{prop:keypart} $\comp{a}_{\lambda} \subseteq \comp{b}$ if and only if $\comp{a}_{\lambda} \preceq \comp{b}$. In particular, $\lambda \subseteq \mu$ if and only if $\comp{a}_{\lambda} {\preceq} \comp{a}_{\mu}$.
\end{proof}

We have the following relations for intersections and unions of diagrams.

\begin{lemma}
  Given weak compositions $\comp{a},\comp{b},\comp{c}$ of length $n$,
  \begin{enumerate}
  \item if $\comp{a} \preceq \comp{b}$ and $\comp{a} \preceq \comp{c}$, then $\comp{a} \preceq \comp{b} \cap \comp{c}$;
  \item if $\comp{a} \preceq \comp{c}$ and $\comp{b} \preceq \comp{c}$, then $\comp{a} \cup \comp{b} \preceq \comp{c}$.
  \end{enumerate}
  \label{lem:cupcap}
\end{lemma}

\begin{proof}
  Suppose $\comp{a} \preceq \comp{b}$ and $\comp{a} \preceq \comp{c}$. By definition of the poset, we have $\comp{a} \subseteq \comp{b}$ and $\comp{a} \subseteq \comp{c}$, whence $\comp{a} \subseteq \comp{b} \cap \comp{c}$. For any cell $x$ in $(\comp{b}\cap \comp{c}) \skews \comp{a}$, if $x$ sits above a cell $y$ of $\comp{a}$, then necessarily the row of $x$ is strictly shorter than that of $y$ in $\comp{b}$ and in $\comp{c}$ {since $\comp{a} \preceq  \comp{b}, \comp{c}$}. Thus, the row of $x$ is strictly shorter than that of $y$ in $\comp{b}\cap \comp{c}$, and so $\comp{a} \preceq \comp{b} \cap \comp{c}$.
  
  Suppose $\comp{a} \preceq \comp{c}$ and $\comp{b} \preceq \comp{c}$. Then $\comp{a} \subseteq \comp{c}$ and $\comp{b} \subseteq \comp{c}$, whence $\comp{a}\cup \comp{b} \subseteq \comp{c}$. For any cell $x$ in $\comp{c} \skews (\comp{a}\cup \comp{b})$, if $x$ sits above a cell $y$ of $\comp{a}\cup \comp{b}$, then necessarily the row of $x$ is strictly shorter than that of $y$ in $\comp{a}$ {or} in $\comp{b}$. Thus, the row of $x$ is strictly shorter than that of $y$ in {$\comp{c}$ since $\comp{a}, \comp{b} \preceq \comp{c}$}, and so $\comp{a} \cup \comp{b} \preceq \comp{c}$.
\end{proof}

In particular, by Lemma~\ref{lem:cupcap}, if $\comp{a},\comp{b} \preceq \comp{a} \cup \comp{b}$, then $\comp{a}\cup \comp{b}$ is the unique least upper bound of $\comp{a}$ and $\comp{b}$, and if $\comp{a}\cap \comp{b} \preceq \comp{a},\comp{b}$, then $\comp{a} \cap \comp{b}$ is the unique greatest lower bound of $\comp{a}$ and $\comp{b}$. However, neither of these conditions needs be the case.

\begin{proposition}
  The key poset is not a lattice.
  \label{prop:no-lattice}
\end{proposition}

\begin{proof}
  Let $\comp{a} = (2,1,1)$ and $\comp{b}=(2,1,2)$. By Lemma~\ref{lem:cupcap}, any greatest lower bound for $\comp{a}$ and $\comp{b}$ must be contained in $\comp{a} \cap \comp{b} = \comp{a}$, but $\comp{a} \not\prec \comp{b}$ by Example~\ref{ex:contain}. Thus any greatest lower bound must {have} rank at most $3$. Let $\comp{c} = (1,1,1)$. Then one can check (visually from Fig.~\ref{fig:no-lattice}) $\comp{c} \prec \comp{a},\comp{b}$, so since $\comp{c}$ has rank $3$, it is a greatest lower bound for $\comp{a}$ and $\comp{b}$. Now let $\comp{d}=(1,0,1)$. Once again {one can check}, $\comp{d} \prec \comp{a},\comp{b}$. However, $\comp{d} \not\prec \comp{c}$ since $c_2 > d_2$ and $d_1 > d_2$ but $c_1 = c_2$. Thus $\comp{c}$ is not unique.

  Dually, Lemma~\ref{lem:cupcap} ensures any least upper bound for $\comp{c}$ and $\comp{d}$ must be contained in $\comp{c} \cup \comp{d} = \comp{c}$, but we have just seen $\comp{d} \not\prec \comp{c}$. Thus any least upper bound has rank at least $4$, and so the earlier check confirms $\comp{a}$ as a least upper bound. Then as $\comp{b}$ is {likewise} also a common upper bound and $\comp{a} \not\prec \comp{b}$, we see $\comp{a}$ is not unique.
\end{proof}

\begin{figure}[ht]
  \begin{displaymath}
    \arraycolsep=2\smcellsize
    \begin{array}{ccccc}
    \vline\smtab{ \ \\ \ \\ \ & \cbox{red}} &
    \vline\smtab{ \ & \cbox{blue} \\ \ \\ \ & \cbox{blue}} &
    \vline\smtab{ \ \\ \cbox{red} \\ \ & \cbox{red}} &
    \vline\smtab{ \ & \cbox{blue} \\ \cbox{blue} \\ \ & \cbox{blue}} &
    \vline\smtab{ \ \\ \cbox{green} \\ \ } \\
    \comp{c} \prec \comp{a} &
    \comp{c} \prec \comp{b} &
    \comp{d} \prec \comp{a} &
    \comp{d} \prec \comp{b} &
    \comp{d} \not\prec \comp{c}
    \end{array}
  \end{displaymath}
  \caption{\label{fig:no-lattice}An example showing the key poset is not a lattice. Here ${\cbox{white}}\in \comp{c}$ or $\comp{d}$, ${\cbox{red}} \in \comp{a}\skews \comp{c}$ or $\comp{a}\skews \comp{d}$, ${\cbox{blue}} \in \comp{b}\skews \comp{c}$ or $\comp{b}\skews \comp{d}$, and ${\cbox{green}} \in \comp{c}\skews \comp{d}$.}
\end{figure}

The set-theoretic unions and intersections are natural candidates for least upper bounds and greatest lower bounds, though as Proposition~\ref{prop:no-lattice} shows, they are not necessarily above or below their constituent parts in the key poset. While this problem cannot be overcome completely, there are derived diagrams that do always lie above or below their constituent parts and which are contained in or contain all other least {upper} or greatest {lower} bounds.

\begin{proposition}
  Given weak compositions $\comp{a}$ and $\comp{b}$, let $\comp{a} \nabla \comp{b}$ denote $\comp{a} \cup \comp{b}$ together with all cells $y$ such that the cell immediately left of $y$ sits under a cell of $(\comp{a} \cup \comp{b}) \skews (\comp{a} \cap \comp{b})$ and some cell left of $y$ either lies in $\comp{a}$ with the cell above in $\comp{b} \skews (\comp{a} \cap \comp{b})$ or lies in $\comp{b}$ with the cell above in $\comp{a}\skews (\comp{a} \cap \comp{b})$; see Fig.~\ref{fig:join}.

  Then $\comp{a},\comp{b} \preceq \comp{a} \nabla \comp{b}$, and for any $\comp{c}$ for which $\comp{a}, \comp{b} \preceq \comp{c}$, we have $\comp{a} \nabla \comp{b} \subseteq \comp{c}$.
  \label{prop:join}
\end{proposition}

\begin{figure}[ht]
  \begin{displaymath}
    \vline\smtab{ \ & \ & \cbox{blue} & \cbox{blue} \\ \\ \ & \ & \ & \cbox{purple} & \cbox{purple}}\hspace{-5\smcellsize}\smtab{ \\ \\ & & & + & + }
  \end{displaymath}
  \caption{\label{fig:join}An illustration of the rule for padding $\comp{a} \cup \comp{b}$ up to $\comp{a} \nabla \comp{b}$. Here cells $\cbox{white}$ lie in $\comp{a}$, cells $\cbox{blue}$ lie in $\comp{b}\skews (\comp{a}\cap \comp{b})$, and cells $\cbox{purple}$ must be added to make $\comp{a} \nabla \comp{b}$.}
\end{figure}

\begin{proof}
  By definition, $\comp{a},\comp{b} \subseteq \comp{a}\cup \comp{b} \subseteq \comp{a} \nabla \comp{b}$. Moreover, if some $x\in (\comp{a} \nabla \comp{b}) \skews \comp{a}$ (resp. $x\in (\comp{a} \nabla \comp{b})\skews \comp{b}$) lies above some $y\in \comp{a}$ (resp. $y\in \comp{b}$), then by construction of $\comp{a} \nabla \comp{b}$, the lower row is strictly longer, so $\comp{a},\comp{b} \preceq \comp{a} \nabla \comp{b}$.

  Let $\comp{c}$ be some weak composition for which $\comp{a},\comp{b} \preceq \comp{c}$. By Lemma~\ref{lem:cupcap}, we have $\comp{a} \cup \comp{b} \preceq \comp{c}$, so, in particular, $\comp{a} \cup \comp{b} \subseteq \comp{c}$. By construction of $\comp{a} \nabla \comp{b}$, any $z \in (\comp{a}\nabla \comp{b})\skews (\comp{a} \cup \comp{b})$ lies in a row $i$ for which there is a higher row $j>i$ that is {strictly} shorter in $\comp{a}$ {and} in $\comp{b}$. Since $\comp{a} \cup \comp{b} \preceq \comp{c}$, we must have $z \in \comp{c}$. In particular, $\comp{a} \nabla \comp{b} \subseteq \comp{c}$.
\end{proof}

{\begin{example} We now interpret Fig.~\ref{fig:join} as a concrete example: Let $\comp{a} = (3,0,2)$ and $\comp{b}=(2,0,4)$, so $\comp{a}\cap\comp{b}=(2,0,2)$, $\comp{a}\cup\comp{b}=(3,0,4)$ and $\comp{a}\nabla\comp{b}=(5,0,4)$.
\label{ex:join}
\end{example}}

\begin{proposition}
  Given weak compositions $\comp{a}$ and $\comp{b}$, let $\comp{a} \Delta \comp{b}$ denote the subset of cells $x \in \comp{a} \cap \comp{b}$ such that for any $y \in \comp{a} \skews (\comp{a}\cap \comp{b})$ and any $z\in \comp{b} \skews (\comp{a}\cap \comp{b})$ lying above $x$, the row of $x$ is strictly longer than that of $y$ in $\comp{a}$ and of $z$ in $\comp{b}$; see Fig.~\ref{fig:meet}.

  Then $\comp{a} \Delta \comp{b} \preceq \comp{a}, \comp{b}$, and for any $\comp{c}$ for which $\comp{c} \preceq \comp{a}, \comp{b}$, we have $\comp{c} \subseteq \comp{a} \Delta \comp{b}$.
  \label{prop:meet}
\end{proposition}

\begin{figure}[ht]
  \begin{displaymath}
    \vline\smtab{ \ & \ & \cbox{blue} & \cbox{blue} \\ \\ \ & \ & \times}
  \end{displaymath}
  \caption{\label{fig:meet}An illustration of the rule for culling $\comp{a} \cap \comp{b}$ down to $\comp{a} \Delta \comp{b}$. Here cells $\cbox{white}$ lie in $\comp{a} \cap \comp{b}$, cells $\cbox{blue}$ lie in $(\comp{a} \cup \comp{b})\skews (\comp{a}\cap \comp{b})$, and the cell marked $\times$ is deleted to make $\comp{a} \Delta \comp{b}$.}
\end{figure}

\begin{proof}
  By definition, $\comp{a} \Delta \comp{b} \subseteq \comp{a} \cap \comp{b} \subseteq \comp{a},\comp{b}$. Moreover, if some ${y}\in \comp{a}\skews (\comp{a} \Delta \comp{b})$ (resp. ${z}\in \comp{b}\skews (\comp{a} \Delta \comp{b})$) lies above some ${x}\in \comp{a} \Delta \comp{b}$, then by construction of $\comp{a} \Delta \comp{b}$, the lower row is strictly longer, so $\comp{a} \Delta \comp{b} \preceq \comp{a},\comp{b}$.

  Let $\comp{c}$ be some weak composition for which $\comp{c} \preceq \comp{a},\comp{b}$. By Lemma~\ref{lem:cupcap}, we have $\comp{c} \preceq \comp{a} \cap \comp{b}$, so, in particular, $\comp{c}\subseteq \comp{a} \cap \comp{b}$. By construction of $\comp{a} \Delta \comp{b}$, any ${w} \in (\comp{a}\cap \comp{b})\skews (\comp{a} \Delta \comp{b})$ must have either a cell ${y}\in \comp{a}\skews (\comp{a} \Delta \comp{b})$ or ${z}\in \comp{b}\skews (\comp{a} \Delta \comp{b})$ above it for which the row of ${y}$ in $\comp{a}$ or the row of ${z}$ in $\comp{b}$ is weakly longer than that of ${w}$, and so we cannot have ${w} \in \comp{c}$. In particular, $\comp{c} \subseteq \comp{a} \Delta \comp{b}$.
\end{proof}

{\begin{example} We now interpret Fig.~\ref{fig:meet} as a concrete example: Let $\comp{a} = (3,0,2)$ and $\comp{b}=(3,0,4)$, so $\comp{a}\cap\comp{b}=(3,0,2)$, $\comp{a}\cup\comp{b}=(3,0,4)$ and $\comp{a}\Delta\comp{b}=(2,0,2)$.
\label{ex:meet}
\end{example}}

%
\section{Polynomials}
%
\label{sec:poly}

Based on the {quasi-Yamanouchi tableaux} of Assaf and Searles \cite{AS18}, Assaf defined {standard key tableaux}  \cite[Definition~3.10]{Ass-W}. We generalize this naturally from key diagrams to skew key diagrams as follows.

\begin{definition}[\cite{Ass-W}]
  A \newword{standard {(skew)} key tableau} is a bijective filling of a {(skew)} key diagram with $1,2,\ldots,n$ such that {rows decrease left to right} and if some entry $i$ is above and in the same column as an entry $k$ with $i<k$, then there is an entry immediately right of $k$, say $j$, and $i<j$. 
  \label{def:SKT}
\end{definition}

\begin{remark}
  There is a natural semistandard analog of Definition~\ref{def:SKT} present in \cite{Ass18} where it is shown that these objects coincide exactly with Mason's semiskyline augmented fillings \cite{Mas09} defined by a triple rule arising from nonsymmetric Macdonald polynomials.
\end{remark}

We denote the set of standard key tableaux of \newword{shape} $\comp{a}$ by $\SKT(\comp{a})$ and extend this naturally to standard skew key tableaux. 

A \newword{saturated chain} in a poset $\mathcal{P}$ is a sequence of elements $p_0 \cover p_1 \cover \cdots \cover p_{n}$.  We may identify standard Young tableaux of shape $\lambda$ with saturated chains in Young's lattice from $\varnothing$ to $\lambda$.

Parallel to the case for Young's lattice, saturated chains from $\varnothing$ to $\comp{a}$ in the key poset precisely correspond to standard key tableaux of shape $\comp{a}$.
Look ahead to Theorem~\ref{thm:skew-chains} for the skew analog.

\begin{theorem}
  For $\varnothing = \comp{a}^{(0)} \cover \comp{a}^{(1)} \cover \cdots \cover \comp{a}^{(n)}{=\comp{a}}$ a saturated chain in the key poset, the standard filling of the key diagram for $\comp{a}$ defined by placing $n-i+1$ into the unique cell of $\comp{a}^{(i)}\skews \comp{a}^{(i-1)}$ is a standard key tableaux.

  Conversely, given $T \in \SKT(\comp{a})$ {with $n$ cells}, setting $\comp{a}^{(0)}=\varnothing$ and, for $i=1,\ldots,n$, setting $\comp{a}^{(i)}$ to be the diagram containing cells labeled $n,n-1,\ldots,n-i+1$ results in a saturated chain from $\varnothing$ to $\comp{a}$ in the key poset. 
  \label{thm:chains}
\end{theorem}

\begin{proof}
  Suppose $\varnothing = \comp{a}^{(0)} \cover \comp{a}^{(1)} \cover \cdots \cover \comp{a}^{(n)}$ {is} a saturated chain in the key poset. Then by Theorem~\ref{thm:cover}, $\comp{a}^{(i-1)}\subset \comp{a}^{(i)}$ and $\comp{a}^{(i)}\skews \comp{a}^{(i-1)}$ has one cell. Thus we may indeed set $\comp{a}$ to be the bijective filling of $\comp{a}$ with $n-i+1$ into the unique cell of $\comp{a}^{(i)}\skews \comp{a}^{(i-1)}$, for $i=1,\ldots,n$. The cover relations ensure cells are added only to the right end of a row, ensuring row entries decrease from left to right. Suppose $i<k$ are in the same column with $i$ above $k$. Restricting our attention to $\comp{a}^{(n-i+1)}$, the last cell added corresponds to the cell with entry $i$. By Theorem~\ref{thm:cover}, no row below that containing entry $i$ can have the same length. In particular, the row of $k$ must be strictly longer. Thus there exists some entry $j<k$ immediately right of $k$, and since $j$ was not the most recently added cell, we have $j> i$. Therefore $\comp{a}$ is a standard key tableau.

  Suppose now $T$ is a standard key tableau of size $n$. Since entries decrease from left to right, the {shape of the} restriction of $T$ to entries $n,n-1,\ldots,n-i+1$ must be a key diagram. Therefore we may define a nested sequence of weak compositions $\varnothing = \comp{a}^{(0)} \subset \comp{a}^{(1)} \subset \cdots \subset \comp{a}^{(n)}$ by setting $\comp{a}^{(i)}$ to be the diagram containing cells of $T$ labeled $n,n-1,\ldots,n-i+1$. Note $T$ restricted to these entries still satisfies the row and column conditions for a standard key tableau. To see that $\comp{a}^{(i-1)}\cover \comp{a}^{(i)}$, the column condition for key tableaux ensures  the cell of $\comp{a}^{(i)}$ with smallest entry, which necessarily corresponds to $\comp{a}^{(i)}\skews \comp{a}^{(i-1)}$, cannot lie above another cell that ends {its} row, since that cell necessarily has {a} larger entry. Thus by Theorem~\ref{thm:cover}, the sequence is a saturated chain {in the key poset}. 
\end{proof}

\begin{figure}[ht]
  \begin{displaymath}
    \arraycolsep=2pt
    \begin{array}{ccccccc}
    \varnothing & \cover & \vline\smtab{ 3 \\ & \\ \\\hline } & \cover & \vline\smtab{ 3 \\ 2 & \\ \\\hline } & \cover & \vline\smtab{ 3 \\ 2 & 1 \\ \\\hline } 
    \end{array}\hspace{6\smcellsize}
    \begin{array}{ccccccc}
    \varnothing & \cover & \vline\smtab{ \\ 3 & \\ \\\hline } & \cover & \vline\smtab{ \\ 3 & 2 \\ \\\hline } & \cover & \vline\smtab{ 1 \\ 3 & 2\\ \\\hline } 
    \end{array}
  \end{displaymath}
  \caption{\label{fig:SKT}Constructing the two standard key tableaux of shape $(0,2,1)$ as saturated chains in the key poset from $\varnothing$ to $(0,2,1)$.}
\end{figure}

For example, Fig.~\ref{fig:SKT} shows the two saturated chains in the key poset from $\varnothing$ to $(0,2,1)$.

The {key polynomials}, indexed by weak compositions, form an important basis for the full polynomial ring. Key polynomials arise as characters of Demazure modules \cite{Dem74a} for the general linear group and coincide with Schubert polynomials \cite{LS82} in the vexillary case \cite{LS90}. Key polynomials are nonsymmetric generalizations of Schur functions, studied combinatorially by Reiner and Shimozono \cite{RS95} and later by Mason \cite{Mas09}, though our perspective follows that of Assaf and Searles \cite{AS18} and Assaf \cite{Ass-W} who define them as the {fundamental slide generating polynomial} for standard key tableaux.

Assaf and Searles introduced the  {fundamental slide polynomials} \cite{AS17}, indexed by weak compositions, that form a basis for the full polynomial ring.
\svw{Given strong compositions $\alpha,\beta$, we say $\beta$ \newword{refines} $\alpha$ if there exist indices $i_1<{\cdots}<i_k$, where $k$ is the length of $\alpha$, such that  for all $1 \le j \le k$
 \begin{displaymath}
  \beta_1 + \cdots + \beta_{i_j} = \alpha_1 + \cdots + \alpha_j.
\end{displaymath}
 For example, $(1,2,2)$ refines $(3,2)$ but does not refine $(2,3)$.}

\begin{definition}[\cite{AS17}]
  For a weak composition $\comp{a}$, the \newword{fundamental slide polynomial} $\fund_{\comp{a}}$ is 
  \begin{equation}
    \fund_{\comp{a}} = \sum_{\substack{\mathrm{flat}(\comp{b}) \ \mathrm{refines} \ \mathrm{flat}(\comp{a}) \\ b_1+\cdots+b_k \geq a_1+\cdots+a_k \ \forall k}} x_1^{b_1} x_2^{b_2} \cdots x_n^{b_n},
    \label{e:slide}
  \end{equation}
  where the sum is over weak compositions $\comp{b}$ such that $\mathrm{flat}(\comp{b})$ refines $\mathrm{flat}(\comp{a})$ and $\comp{b}$ dominates $\comp{a}$.
  \label{def:slide}
\end{definition}

For a standard key tableau $T$, say  $i$ is a \newword{descent} of $T$ if $i+1$ lies weakly \emph{right} of $i$. We assign a \newword{weak descent composition} for $T$, defined in \cite[Definition~3.12]{Ass-W}, that will index the corresponding fundamental slide polynomial. 

\begin{definition}[\cite{Ass-W}]
  For a standard key tableau $T$, let $(\tau^{(k)}|{\cdots}|\tau^{(1)})$ be the partitioning of the decreasing word $n \cdots 2 1$ broken between $i+1$ and $i$ precisely whenever $i$ is a descent of $T$. Set $t^{\prime}_i$ to be the lowest row index in $T$ of a letter in $\tau^{(i)}$. Set $t_k = t^{\prime}_k$ and, for $i<k$, set $t_i = \min(t^{\prime}_i , t_{i+1}-1)$. Define the \newword{weak descent composition of $T$}, denoted by $\des(T)$, by $\des(T)_{t_i} = |\tau^{(i)}|$ and all other parts are zero if all $t_i>0$; and $\des(T) = \varnothing$ otherwise.
  \label{def:des}
\end{definition}


\begin{example}
  The left tableau in Fig.~\ref{fig:SKT} has a descent only at $2$, so the partitioning is $(\tau^{(2)} | \tau^{(1)}) = (3 | 21)$. Thus $t_2 = t^{\prime}_2 = 3$, the row index of $3$, $t^{\prime}_1 = 2$, the row index of $2$ and $1$, and $t_1 = \min(t^{\prime}_1,t_2-1) = \min(2,2) = 2$. These determine the \emph{positions} of the descents, and so the weak descent composition is $(0,|\tau^{(1)}|,|\tau^{(2)}|) = (0,2,1)$.

  The right tableau in Fig.~\ref{fig:SKT} has a descent only at $1$, so the partitioning is now $(\tau^{(2)} | \tau^{(1)}) = (32 | 1)$. Thus $t_2 = t^{\prime}_2 = 2$, the row index of $2$ and $3$, {$t^{\prime}_1 = 3$}, the row index of $1$, and $t_1 = \min(t^{\prime}_1,t_2-1) = \min({3,}1) = 1$. Thus the weak descent composition is $(|\tau^{(1)}|,|\tau^{(2)}|,0) = (1,2,0)$.
  \label{ex:des}
\end{example}

We take \cite[Corollary~3.16]{Ass-W} as our definition for key polynomials.

\begin{definition}[\cite{Ass-W}]
  For a weak composition $\comp{a}$, the \newword{key polynomial} $\key_{\comp{a}}$ is 
  \begin{equation}
    \key_{\comp{a}} = \sum_{T \in \SKT(\comp{a})} \fund_{\des(T)},
    \label{e:key}
  \end{equation}
  where the sum is over all standard key tableaux of shape $\comp{a}$ for which the weak descent composition is not $\varnothing$.
  \label{def:key}  
\end{definition}

For example, from Example~\ref{ex:des}, we have
\[ \kappa_{(0,2,1)} = \fund_{(0,2,1)} + \fund_{(1,2,0)} . \]

Note working with semistandard in place of standard objects would replace the fundmental slide polynomials with their monomial expansions, giving a far less compact notation.

Composing the bijective correspondence defined by Assaf and Searles between quasi-Yamanouchi Kohnert tableaux and quasi-Yamanouchi Young tableaux \cite[Theorem~4.6]{AS18} with the bijective correspondence defined by Assaf \cite[Theorem~3.15]{Ass-W} between those and the corresponding standard tableaux yields
a bijective proof that for $\comp{a}=(0,\ldots,0,\lambda_{\ell},\ldots,\lambda_1)$ {of} length $n$, we have
\begin{equation}
  \key_{(0,\ldots,0,\lambda_{\ell},\ldots,\lambda_1)} = s_{\lambda}(x_1,x_2,\ldots,x_n).
\end{equation}
Thus key polynomials generalize Schur polynomials. Moreover, key polynomials \newword{stabilize} to Schur functions as 
\begin{equation}
  \lim_{m \rightarrow\infty} \key_{0^m\times a}(x_1,\ldots,x_m,0,\ldots) = s_{\lambda}(X),
\end{equation}
where $0^m\times \comp{a}$ is $\comp{a}$ with $m$ zeros prepended and $\lambda$ is the partition reordering of $\comp{a}$.


%
\section{Littlewood--Richardson rules}
%
\label{sec:skew}

Assaf considers {standard skew key tableaux} in \cite[Definition~4.7]{Ass-W}, defined for any pair of weak compositions $\comp{a} \subset {\comp{d}}$. However, the positivity result for the corresponding {skew key polynomials} \cite[Theorem~4.10]{Ass-W}, 
is only proved for the limited case when the smaller weak composition is a \emph{partition}. Examples show  these are not the only cases where nonnegativity holds, with many overlooked examples arising naturally from geometric contexts. As we shall see, the key to positivity lies in the key poset. We begin by generalizing Theorem~\ref{thm:chains} to skew {key} diagrams.

\begin{theorem}
  For $\comp{a} = \comp{a}^{(0)} \cover \comp{a}^{(1)} \cover \cdots \cover \comp{a}^{(n)} = \comp{d}$ a saturated chain in the key poset, the standard filling of the skew key diagram for $\comp{d}\skews \comp{a}$ defined by placing $n-i+1$ into the unique cell of $\comp{a}^{(i)}\skews \comp{a}^{(i-1)}$ is a standard skew key tableaux.

  Conversely, for $T \in \SKT(\comp{d}\skews \comp{a})$ with $n$ cells, setting $\comp{a}^{(0)}={\comp{a}}$ and, for $i=1,\ldots,n$, setting $\comp{a}^{(i)}$ to be the diagram with cells of $\comp{a}$ skewed and containing cells labeled $n,n-1,\ldots,n-i+1$ gives a saturated chain from $\comp{a}$ to $\comp{d}$ in the key poset. 
  \label{thm:skew-chains}
\end{theorem}

\begin{proof}
  Both the poset cover relations and the key tableaux row and column conditions are local, so this follows from Theorem~\ref{thm:chains} by restricting attention to shapes from $\comp{a}$ onward and to cells with entries $1,2,\ldots,\rk(\comp{d})-\rk(\comp{a})$.
\end{proof}

Following \cite{Ass-W}, extend Definition~\ref{def:des} directly to \svw{standard skew} key tableaux
of shape $\comp{d} \skews \comp{a}$ by regarding the cells of $\comp{a}$ as labeled $n,n-1,\ldots,n-r+1$, where $n = |\comp{d}|$ and $r = |\comp{a}|$.

\begin{definition}
  For weak compositions $\comp{a} \prec \comp{d}$, the \newword{skew key polynomial} $\key_{\comp{d}\skews \comp{a}}$ is 
  \begin{equation}
    \key_{\comp{d} \skews \comp{a}} = \sum_{T \in {\SKT(\comp{d}\skews \comp{a})}} \fund_{\des(T)},
    \label{e:skewkey}
  \end{equation}
  where the sum is over all standard \svw{skew} key tableaux of shape \svw{$\comp{d}\skews \comp{a}$} for which the weak descent composition is not $\varnothing$.
  \label{def:skewkey}
\end{definition}

Note  unlike \cite[Definition~4.8]{Ass-W}, we define skew key polynomials only for comparable elements of the key poset. However, the special case of skewing by an increasing composition, for which the positivity in \cite[Theorem~4.10]{Ass-W} holds, conforms with this more restrictive definition by Proposition~\ref{prop:keypart}. 

Since key polynomials are a basis for all polynomials,  {define} \newword{weak composition Littlewood--Richardson coefficients} $c_{\comp{a},\comp{b}}^{\comp{d}}$ as the key expansion of skew key polynomials,
\begin{equation}
  \key_{\comp{d} \skews \comp{a}} = \sum_{\comp{b}} c_{\comp{a},\comp{b}}^{\comp{d}} \key_{\comp{b}} .
\end{equation}
A priori, these coefficients are \emph{integers}. In fact, we will show that they are \emph{nonnegative} integers, and so skew key polynomials are \newword{key positive}.

To prove nonnegativity of the {weak} composition Littlewood--Richardson coefficients, we utilize {weak dual equivalence} \cite{Ass-W}, a polynomial generalization of dual equivalence that consolidates standard key tableaux into equivalence classes, each of which corresponds to a single key polynomial. 

Extending earlier notation, given a weak composition $\comp{a}$ of rank $n$ and integers $1\leq h \leq i\leq n$, let $\comp{a}_{(h,i)}$ be the weak composition obtained by deleting the first $h-1$ and last $n-i$ pieces from $\comp{a}$.

\begin{example}
  Let $\comp{a} = (0,3,2,0,3,1)$. Then $\comp{a}_{(3,7)} = (0,1,2,0,2,0)$, corresponding to {deleting} the first $2$ pieces both of which come from $a_2$ and last $2$ pieces one of which comes from $a_6$ and the remaining from $a_5$.
\end{example}

\begin{definition}[\cite{Ass-W}]
  Let $\mathcal{A}$ be a finite set, and let $\des$ be a map from $\mathcal{A}$ to weak compositions of rank $n$. A \newword{weak dual equivalence for $(\mathcal{A},\des)$} is a family of involutions $\{\psi_i\}_{1<i<n}$ on $\mathcal{A}$ such that
  \renewcommand{\theenumi}{\roman{enumi}}
  \begin{enumerate}
  \item For all $i-h \leq 3$ and $T \in \mathcal{A}$, there exists a weak composition $\comp{a}$ such that
    \[ \sum_{U \in [T]_{(h,i)}} \fund_{\des_{(h-1,i+1)}(U)} = \key_{\comp{a}}, \]
    where $[T]_{(h,i)}$ is the equivalence class generated by $\psi_h,\ldots,\psi_i$.
    
  \item For all $|i-j| \geq 3$ and all $T \in\mathcal{A}$, we have $\psi_{j} \psi_{i}(T) = \psi_{i} \psi_{j}(T)$.
  \end{enumerate}
  \label{def:wdeg}
\end{definition}

\svw{Define the \newword{simple involutions} $s_i$ for $i=1,2,\ldots,n-1$ on standard fillings of rank $n$ that interchange $i$ and $i+1$. }
Define the \newword{braid involutions} $b_i$ for $i=2,\ldots,n-1$ on standard fillings of rank $n$ for which exactly one of $i-1$ or $i+1$ lies in the same row as $i$ by cycling entries $i-1,i,i+1$ in the unique way that maintains this condition. For example, $b_2$ will exchange the two standard key tableaux of shape $(0,2,1)$ shown in Fig.~\ref{fig:SKT}. We combine the braid involutions with \svw{the simple involutions} to give the elementary weak dual equivalence involutions in \cite[Definition~3.21]{Ass-W}.

\begin{definition}[\cite{Ass-W}]
  Define \newword{elementary weak dual equivalence involutions}, denoted by $d_i$, on  \svw{standard skew} key tableaux that act by
  \begin{equation}
    d_i(T) = \left\{ \begin{array}{rl}
      b_{i} \cdot T & \text{if exactly one of $i-1,i+1$ lies in the row of $i$, } \\
      s_{i-1} \cdot T & \text{else if $i+1$ lies between $i$ and $i-1$ in reading order, } \\
      s_{i} \cdot T & \text{else if $i-1$ lies between $i$ and $i+1$ in reading order, } \\
      T & \text{otherwise} 
    \end{array} \right.
    \label{e:wede}
  \end{equation}
  where we take {the reading order to be taking the entries in each \emph{column},} bottom to top and left to right.
  \label{def:wede}
\end{definition}

For examples of the elementary weak \svw{dual equivalence} involutions, see Fig.~\ref{fig:wdeg}.

\begin{figure}[ht]  
  \begin{center}
    \begin{tikzpicture}[xscale=1.5,yscale=1]
      \node at (0,0)   (a) {$\smtab{5 & 2 \\ 4 \\ 3 & 1 }$}; 
      \node at (1,0)   (b) {$\smtab{5 & 3 \\ 4 \\ 2 & 1 }$}; 
      \node at (2,0)   (c) {$\smtab{5 & 4 \\ 3 \\ 2 & 1 }$}; 
      \node at (3,0)   (d) {$\smtab{5 & 4 \\ 1 \\ 3 & 2 }$}; 
      \node at (4,0)   (e) {$\smtab{5 & 3 \\ 1 \\ 4 & 2 }$}; 
      \draw[thick,color=blue  ,<->] (a.008) -- (b.172) node[midway,above] {$d_2$} ;
      \draw[thick,color=violet,<->] (a.352) -- (b.188) node[midway,below] {$d_3$} ;
      \draw[thick,color=purple,<->] (b) -- (c) node[midway,above] {$d_4$} ;
      \draw[thick,color=blue  ,<->] (c) -- (d) node[midway,above] {$d_2$} ;
      \draw[thick,color=violet,<->] (d.008) -- (e.172) node[midway,above] {$d_3$} ;
      \draw[thick,color=purple,<->] (d.352) -- (e.188) node[midway,below] {$d_4$} ;
    \end{tikzpicture}
  \end{center}
  \caption{\label{fig:wdeg}The elementary \svw{weak} dual equivalence involutions on $\SKT(2,1,2)$.}
\end{figure}

Assaf \cite[Theorem~3.25]{Ass-W} showed  these are well-defined involutions on standard key tableaux, that all standard key tableaux of fixed shape fall into a single equivalence class, and that this gives an example of a weak dual equivalence.

\begin{theorem}[\cite{Ass-W}]
  Given a weak composition $\comp{a}$, the involutions $d_i$ give a weak dual equivalence for $(\SKT(\comp{a}),\des)$ consisting of a single equivalence class. 
  \label{thm:deg-key}
\end{theorem}

Moreover, under certain stability conditions, the converse holds. That is, by \cite[Theorem~3.29]{Ass-W}, any weak dual equivalence is essentially this and, on the level of generating polynomials, we have the following.

\begin{theorem}[\cite{Ass-W}]
  If there exists a weak dual equivalence for $(\mathcal{A},\mathrm{des})$ for which $\des(T)\neq\varnothing$ for every $T\in\mathcal{A}$, then
  \begin{displaymath}
    \sum_{T \in \mathcal{A}} \fund_{\mathrm{des}(T)}
  \end{displaymath}
  is {key positive}.
  \label{thm:wdeg}
\end{theorem}

The condition that the weak descent composition is nonempty for every element can often be circumvented if the polynomials under consideration stabilize. Assaf uses this along with the elementary weak dual equivalence involutions to prove the special case of the following when $\comp{a}$ is weakly increasing \cite[Theorem~4.10]{Ass-W}, which follows from Theorem~\ref{thm:skewpos} by Proposition~\ref{prop:keypart}.

\begin{theorem}
  For $\comp{a} \prec \comp{d}$ in the key poset, the {weak} composition Littlewood--Richardson coefficient $c_{\comp{a},\comp{b}}^{\comp{d}}$ is the number of weak dual equivalence classes of $\SKT(\comp{d}\skews$ $\comp{a})$ isomorphic to $\SKT(\comp{b})$. In particular, skew key polynomials are key positive.
  \label{thm:skewpos}
\end{theorem}

\begin{proof}
  By Theorem~\ref{thm:deg-key}, the elementary weak dual equivalence involutions of Definition~\ref{def:wede} give a weak dual equivalence on $\SKT(\comp{d})$. For $\comp{a} \prec \comp{d}$ in the key poset, by Theorem~\ref{thm:skew-chains}, we have $\SKT(\comp{d}\skews\comp{a}) \subset \SKT(\comp{d})$. Since Definition~\ref{def:wdeg} is completely local, these same involutions {restrict} to a weak dual equivalence on $\SKT(\comp{d}\skews\comp{a})$. 

  By Definition~\ref{def:des}, it is clear there exists some nonnegative integer $m$ for which $\des(T)\neq\varnothing$ for every $T\in\SKT(0^m\times\comp{d})$ (the minimal such $m$ is given explicitly in \cite[Theorem~A.6]{AS18}). Therefore $\SKT(0^m\times(\comp{d}\skews\comp{a})) \subset \SKT(0^m\times\comp{d})$ also satisfies the hypotheses of Theorem~\ref{thm:wdeg}, and so we conclude $\key_{0^m\times(\comp{d}\skews\comp{a})}$ is a nonnegative sum of key polynomials. Setting the first $m$ variables to $0$, {we have that} $\key_{\comp{d}\skews\comp{a}}$ is also key positive.
\end{proof}

In fact, we can use the poset structure to prove  this result is tight.

\begin{theorem}
  For $\comp{a} \subset \comp{d}$ such that $\comp{a} \not\prec \comp{d}$, there exists a weak composition $\comp{b}$ for which the {weak} composition Littlewood--Richardson coefficient $c_{\comp{a},\comp{b}}^{\comp{d}}$ is negative. 
  \label{thm:skewneg}
\end{theorem}

\begin{proof}
  Suppose $a_i \leq d_i$ for all $i$ but $\comp{a} \not\prec \comp{d}$. Let $i$ be the largest index for which there exists an index $j>i$ such that $a_i > a_j$ but, contrary to Definition~\ref{def:above}, $d_j \leq d_i$, and take $j$ largest among these. First consider the minimal case when both $\comp{a}$ and $\comp{d}$ have at most $2$ nonzero parts. Generically, we may assume $\comp{a} = (0^{i-1},q,0^{j-i-1},p)$ for $q>p \geq 0$ and $\comp{d} = (0^{i-1},m,0^{j-i-1},n)$ where $n \geq m$ and $n>p$, as illustrated in Fig.~\ref{fig:negative}. In this case, {using Definition~\ref{def:skewkey},} we have the expansion
  \[ \key_{\comp{d}\skews\comp{a}} = \key_{(0^{i-1},m-q,0^{j-i-1},n-p)} - \key_{(0^{i-1},n-p,0^{j-i-1},m-q)} + \text{lower terms}, \]
  where the lower terms do not use the variable $x_{j}$. Since $n-p > m-q$, the two terms with $x_j$ do not cancel, and so the key expansion is not nonnegative.

  In the general case, the term $\key_{\hat{\comp{b}}}$ appears with coefficient $1$ where $\hat{b}_k = d_k - a_k$, and the term $\key_{\comp{b}}$ appears with coefficient $-1$ where $b_k = \hat{b}_k$ for $k \neq i,j$ and $b_i = \hat{b}_j$ and $b_j = \hat{b}_i$. Since $b_i > b_j$, these terms do not cancel, and so $c_{\comp{a},\comp{b}}^{\comp{d}}$ is negative. 
\end{proof}

  \begin{figure}[ht]
    \begin{displaymath}
      \vline\smtab{ %
        \cbox{white} & \cdots & \cbox{white} & \cbox{blue} & \cbox{blue} & \cbox{blue} & \cbox{blue} & \cbox{blue} \\ \\
        \cbox{white} & \cdots & \cbox{white} & \cbox{white} & \cbox{blue} & \cbox{blue} & \cbox{blue} \\ &}
      \hspace{-8\smcellsize}\smtab{ %
         &  &  &  &  & \cdots &  & \cdots \\ \\
         &  &  &  &  & \cdots &  \\ &}
    \end{displaymath}
    \caption{\label{fig:negative}A minimal instance of $\comp{a} \subset \comp{d}$ but $\comp{a} \not\prec \comp{d}$ in the key poset. Here, cells $\cbox{white}$ lie in $\comp{a} \subseteq \comp{d}$ and cells $\cbox{blue}$ lie in $\comp{d} \skews \comp{a}$.}
  \end{figure}

Taken together, Theorems~\ref{thm:skewpos} and \ref{thm:skewneg} show  the key poset precisely characterizes the skew key polynomials with nonnegative key polynomial expansion.

The {flagged Schur polynomials} are another polynomial generalization of Schur functions originally defined by Lascoux and Sch\"{u}tzenberger~\cite{LS82} and studied further by Wachs~\cite{Wac82}. Given partitions $\mu\subset\lambda$ and a {\newword{flag}} ${\bf b}=(b_1\le b_2\le {\cdots \le} b_{\ell_{\lambda}})$, the \newword{flagged skew Schur polynomial} $S_{\lambda\skews\mu,{\bf b}}$ is the sum of the monomials corresponding to semistandard Young tableaux of skew shape $\lambda\skews\mu$ with entries in row $i$ bounded above by $b_i$. Taking $\mu=\varnothing$ gives the (nonskew) flagged Schur polynomial.

In addition to the tableaux definition, Assaf and Bergeron \cite{AB} describe flagged (skew) Schur functions in terms of \emph{flagged $(\mathcal{P},\rho)$-partitions}, a description that immediately gives the fundamental slide expansion as well.

Reiner and Shimozono \cite[Theorem~23]{RS95} show  every flagged Schur polynomial $S_{\lambda,{\bf b}}$ is equal to a single key polynomial $\key_{\comp{a}}$ for some weak composition $\comp{a}$, though the converse does not hold. That is, key polynomials are more general than flagged Schur polynomials.

Moreover, Reiner and Shimozono \cite[Theorem~20]{RS95} give a \newword{flagged Littlewood--Richardson rule} showing  any flagged skew Schur polynomial is a positive sum of key polynomials, that is
\begin{equation}
  S_{\lambda\skews\mu,{\bf b}} = \sum_{\comp{a}} c_{\comp{a},\mu}^{\lambda} \key_{\comp{a}},
\end{equation}
where $c_{\comp{a},\mu}^{\lambda}$ counts the number of $\lambda\skews\mu$-compatible tableaux whose left-nil key with respect to $\lambda-\mu$ is $\comp{a}$; see \cite{RS95} for definitions and details.

By Theorems~\ref{thm:skewpos} and \ref{thm:skewneg}, it follows that each flagged skew Schur polynomial $S_{\lambda\skews\mu,{\bf b}}$ is equal to a skew key polynomial $\key_{\comp{d}\skews\comp{a}}$ for some weak compositions $\comp{a} \prec \comp{d}$, though once again the converse does not hold. That is, skew key polynomials are more general than flagged skew Schur polynomials, and so we obtain a maximal generalization of the flagged Littlewood--Richardson rule.

\section{Acknowledgments}\label{sec:ack} \svw{The authors would like to thank the referees for useful suggestions and careful reading.}
%
%

\bibliographystyle{plain} 

\begin{thebibliography}{10}

\bibitem{Ass18}
Sami Assaf.
\newblock Nonsymmetric {M}acdonald polynomials and a refinement of
  {K}ostka--{F}oulkes polynomials.
\newblock {\em Trans. Amer. Math. Soc.}, 370(12):8777--8796, 2018.

\bibitem{AB}
Sami Assaf and Nantel Bergeron.
\newblock Flagged {$(\mathcal{P},\rho)$}-partitions.
\newblock {\em European J. Combin.}, 86:103085, 17, 2020.

\bibitem{AQ18}
Sami Assaf and Danjoseph Quijada.
\newblock A {P}ieri rule for key polynomials.
\newblock {\em S\'{e}m. Lothar. Combin.}, 80B:Art. 78, 12, 2018.

\bibitem{AS17}
Sami Assaf and Dominic Searles.
\newblock Schubert polynomials, slide polynomials, {S}tanley symmetric
  functions and quasi-{Y}amanouchi pipe dreams.
\newblock {\em Adv. Math.}, 306:89--122, 2017.

\bibitem{AS18}
Sami Assaf and Dominic Searles.
\newblock Kohnert tableaux and a lifting of quasi-{S}chur functions.
\newblock {\em J. Combin. Theory Ser. A}, 156:85--118, 2018.

\bibitem{Ass-W}
Sami~H. Assaf.
\newblock Weak dual equivalence for polynomials.
\newblock {\em Ann. Comb.}, to appear.

\bibitem{BLvW11}
C.~Bessenrodt, K.~Luoto, and S.~van Willigenburg.
\newblock Skew quasisymmetric {S}chur functions and noncommutative {S}chur
  functions.
\newblock {\em Adv. Math.}, 226(5):4492--4532, 2011.

\bibitem{Dem74a}
Michel Demazure.
\newblock D\'esingularisation des vari\'et\'es de {S}chubert
  g\'en\'eralis\'ees.
\newblock {\em Ann. Sci. \'Ecole Norm. Sup. (4)}, 7:53--88, 1974.
\newblock Collection of articles dedicated to Henri Cartan on the occasion of
  his 70th birthday, I.

\bibitem{Dem74}
Michel Demazure.
\newblock Une nouvelle formule des caract\`eres.
\newblock {\em Bull. Sci. Math. (2)}, 98(3):163--172, 1974.

\bibitem{HLMvW11-2}
J.~Haglund, K.~Luoto, S.~Mason, and S.~van Willigenburg.
\newblock Quasisymmetric {S}chur functions.
\newblock {\em J. Combin. Theory Ser. A}, 118(2):463--490, 2011.

\bibitem{HLMvW11}
J.~Haglund, K.~Luoto, S.~Mason, and S.~van Willigenburg.
\newblock Refinements of the {L}ittlewood-{R}ichardson rule.
\newblock {\em Trans. Amer. Math. Soc.}, 363(3):1665--1686, 2011.

\bibitem{Kuono}
T.~Kouno.
\newblock Decomposition of tensor products of {D}emazure crystals.
\newblock {\em J. Algebra}, 546:641--678, 2020.

\bibitem{LS82}
Alain Lascoux and Marcel-Paul Sch{\"u}tzenberger.
\newblock Polyn\^omes de {S}chubert.
\newblock {\em C. R. Acad. Sci. Paris S\'er. I Math.}, 294(13):447--450, 1982.

\bibitem{LS90}
Alain Lascoux and Marcel-Paul Sch{\"u}tzenberger.
\newblock Keys \& standard bases.
\newblock In {\em Invariant theory and tableaux ({M}inneapolis, {MN}, 1988)},
  volume~19 of {\em IMA Vol. Math. Appl.}, pages 125--144. Springer, New York,
  1990.

\bibitem{LR34}
D.~E. Littlewood and A.~R. Richardson.
\newblock Group characters and algebra.
\newblock {\em Philos. Trans. R. Soc. Lond. Ser. A}, 233:99--141, 1934.

\bibitem{LMvW}
K.~Luoto, S.~Mykytiuk, and S.~van Willigenburg.
\newblock {\em An introduction to quasisymmetric Schur functions - Hopf
  algebras, quasisymmetric functions, and Young composition tableaux}.
\newblock Springer, 2013.

\bibitem{Mas09}
Sarah Mason.
\newblock An explicit construction of type {A} {D}emazure atoms.
\newblock {\em J. Algebraic Combin.}, 29(3):295--313, 2009.

\bibitem{RS95}
Victor Reiner and Mark Shimozono.
\newblock Key polynomials and a flagged {L}ittlewood-{R}ichardson rule.
\newblock {\em J. Combin. Theory Ser. A}, 70(1):107--143, 1995.

\bibitem{Wac82}
Michelle~L. Wachs.
\newblock Flagged {S}chur functions, {S}chubert polynomials, and symmetrizing
  operators.
\newblock {\em J. Combin. Theory Ser. A}, 40(2):276--289, 1985.

\end{thebibliography}

\end{document}